\documentclass[preprint,12pt,number]{elsarticle}
\usepackage[T1]{fontenc}
\usepackage{lmodern,microtype}
\usepackage{amsmath,amssymb,amsthm,mathtools,enumitem,booktabs}
\usepackage[hidelinks]{hyperref}
\usepackage{bookmark}
\biboptions{sort&compress}
\journal{Discrete Mathematics}
\allowdisplaybreaks
\newtheorem{theorem}{Theorem}
\newtheorem{claim}{Claim}[section]
\renewcommand{\theclaim}{\arabic{claim}}

\newenvironment{theoremproof}
  {\par%
   \pushQED{\qed}\ignorespaces}
  {\popQED\par}
\newenvironment{claimproof}
  {%
   \begin{proof}[Proof of Claim~\theclaim]}
  {\end{proof}}

\theoremstyle{remark}
\newtheorem*{exampleintro}{Example}
\begin{document}
\begin{frontmatter}
\title{Ramsey numbers of fixed-length cycles, paths and stars in set-colored triple systems}
\author[ccnu]{Xuanyu Cao}
\ead{2025112111@mails.ccnu.edu.cn}
\author[ccnu,hubei]{Chunxiang Wang\corref{cor}}
\ead{wcxiang@mail.ccnu.edu.cn}
\cortext[cor]{Corresponding author.}
\address[ccnu]{School of Mathematics and Statistics, Central China Normal University,
  152 Luoyu Road, Wuhan 430079, China}
\address[hubei]{Hubei Key Laboratory of Mathematical Sciences, Central China Normal University,
  Wuhan 430079, China}

\begin{abstract}
Let $R_3(H;r,s,t)$ be the least order of a complete triple system in which
every assignment of $s$ of $r$ colors to each triple forces a copy of $H$
whose edges share at least $t$ colors. Fix $r\ge s\ge t\ge1$ and let $n$
count edges. We prove that the Ramsey numbers of the linear path and the
linear cycle of fixed length $n$ are both $\gamma n+o(n)$ as $n\to\infty$, with
$\gamma=2$ if and only if $r=s$. For tight paths the corresponding
formula is $\tau n+o(n)$, with the same coefficient for tight cycles
when $3\mid n$. To obtain exactly $n$ edges in a linear cycle, we choose
the parity of the connections and vary the lengths of the paths they join.
For a linear star $S_n$, whose edges meet only at the center,
$R_3(S_n;r,s,t)=cn+O(n^{2/3})$ for fixed colors. For fixed $n\ge2$
and $s,t$ as $r\to\infty$, it is asymptotic to
$\sqrt{6a_n\binom{r}{t}/\binom{s}{t}}$, where $a_n=n(n-1)$ for odd
$n$ and $a_n=n(2n-3)/2$ for even $n$.
For a tight star, whose edges share a pair, the leading term is
$(n-1)\binom{r}{t}/\binom{s}{t}$ in both regimes. For a fixed linear path of
$n\ge3$ edges and fixed $s\ge t\ge2$, the upper bound is
$(2000h/1307+o(1))\binom{r}{t}/\binom{s}{t}$ as $r\to\infty$,
where $h=\lfloor(n-1)/2\rfloor$.
\end{abstract}
\end{frontmatter}
\thispagestyle{plain}

\section{Introduction}\label{sec:intro}
Omidi and Shahsiah~\cite{intro:OS} determined the two-color Ramsey
numbers of linear paths and cycles in triple systems, that is, in
hypergraphs whose edges are triples. Throughout, $n$ counts edges.
The \emph{linear path} $P_n^{\mathrm{lin}}$ has edges
$\{x_{i-1},y_i,x_i\}$ for $1\le i\le n$ on $2n+1$ distinct vertices.
For $n\ge3$, identifying $x_n$ with $x_0$ gives the \emph{linear cycle}
$C_n^{\mathrm{lin}}$. These paths and cycles are also called loose.
For $n\ge3$ the two-color Ramsey number of the path is
$2n+\lfloor(n+1)/2\rfloor$ and that of the cycle is one smaller,
so they have the same leading term.

Conlon, Fox, He, Mubayi, Suk and Verstra\"ete~\cite{intro:CFHMSV}
studied Ramsey numbers when each edge receives a set of colors and
all edges of a copy must share one color. Chen, Ding, Fan, Hu, Jiang
and Wang~\cite{density:CDFHJW} studied the requirement of several
common colors for ordinary graphs, including paths and even cycles.
For triples this raises two questions. Does the equality of the leading
coefficients for paths and cycles persist under set-coloring? Can one
still prescribe the cycle length when consecutive triples overlap in
two vertices? Theorem~\ref{thm:paths-cycles} answers both questions,
with a divisibility condition in the second case.

Write $[N]=\{1,\ldots,N\}$, and let $\binom Xj$ denote the family
of $j$-element subsets of a finite set $X$. The complete triple system
$K_N^{(3)}$ has vertex set $[N]$ and edge set $\binom{[N]}3$.
We also use the term \emph{$3$-graph} for a triple system, write
$v(H)$ and $e(H)$ for its numbers of vertices and edges, and write
$xy$ and $xyz$ for the sets $\{x,y\}$ and $\{x,y,z\}$.
Fix integers $r\ge s\ge t\ge1$. An \emph{$(r,s)$-coloring} assigns
each triple $e$ a label $\chi(e)\in\binom{[r]}s$. A copy of $H$
is an injective image of its vertices carrying each edge to a host
edge. It is \emph{$t$-common} if the intersection of its edge labels
has size at least $t$. Let $R_3(H;r,s,t)$ be the least $N$ such that
every $(r,s)$-coloring of $K_N^{(3)}$ contains a $t$-common copy of $H$.
The ordinary $k$-color Ramsey number is the case $(r,s,t)=(k,1,1)$.

For $T\in\binom{[r]}t$, the \emph{support} $H_T$ is the spanning
$3$-graph formed by the triples $e$ with $T\subseteq\chi(e)$.
A copy is $t$-common precisely when it lies in one of these supports.
There are $\binom{r}{t}$ supports, counted by their color sets, and
each triple belongs to exactly $\binom{s}{t}$ of them. Thus
$\sum_T e(H_T)=\binom{s}{t}\binom{N}{3}$.

In the \emph{tight path} $P_n^{\mathrm{tight}}$ the edges are
$z_i z_{i+1}z_{i+2}$ for $1\le i\le n$ on $n+2$ distinct vertices.
In the \emph{tight cycle} $C_n^{\mathrm{tight}}$ they are the consecutive
triples of a cyclic ordering of $n\ge4$ distinct vertices.

Axenovich, Gy\'arf\'as, Liu and Mubayi~\cite{ext:AGLM} studied the
ordinary multicolor Ramsey numbers of two triples meeting in one
vertex or in a pair. The first theorem treats both configurations
with any number of edges. In the \emph{linear star} $S_n$, the edges
are $xa_ib_i$, with all $2n$ vertices $a_i,b_i$ distinct and different
from $x$. In the \emph{tight star} $B_n$, they are $xyz_i$, with the
$n$ distinct vertices $z_i$ outside the fixed pair $xy$. Thus $S_n$
has $2n+1$ vertices and $B_n$ has $n+2$. We consider both a growing
star and a fixed star with an increasing number of colors.

\begin{theorem}\label{thm:stars}
\leavevmode
\begin{enumerate}[label=(\roman*)]
\item Fix $r\ge s\ge t\ge1$. There is a constant $c=c(r,s,t)$ in
$[2,\sqrt{6\,\binom{r}{t}/\binom{s}{t}}]$ such that, as $n\to\infty$,
\[
 R_3(S_n;r,s,t)=cn+O_{r,s,t}(n^{2/3}),\;
 R_3(B_n;r,s,t)=(n-1)\frac{\binom{r}{t}}{\binom{s}{t}}+O_{r,s,t}(1).
\]
For every $\varepsilon>0$, finitely many rational linear programs give
an interval of length at most $\varepsilon$ containing $c$. The
ordinary two-color coefficient satisfies $7/3\le c(2,1,1)\le3$.
\item Fix $n\ge2$ and $s\ge t\ge1$. Write $a_n=n(n-1)$ for odd $n$ and
$a_n=n(2n-3)/2$ for even $n$. As $r\to\infty$,
\[
 R_3(S_n;r,s,t)\sim\sqrt{6a_n\,\binom{r}{t}/\binom{s}{t}},\qquad
 R_3(B_n;r,s,t)\sim(n-1)\frac{\binom{r}{t}}{\binom{s}{t}}.
\]
\end{enumerate}
\end{theorem}

To obtain a cycle of fixed length, one must join the endpoints of a path and control the number of edges used to close it. Allen, B\"ottcher, Cooley and
Mycroft~\cite{cyextra:ABCM} embedded tight cycles of fixed lengths
divisible by three in dense hypergraphs. In the proof of the next
theorem we use their result in the support of a single $t$-set. For
linear cycles we construct connecting paths of two opposite parities
and vary the lengths of the intervening paths in steps of two.

\begin{theorem}\label{thm:paths-cycles}
Fix $r\ge s\ge t\ge1$. There are constants $\gamma=\gamma_{r,s,t}$
and $\tau=\tau_{r,s,t}$ with $2\le\gamma\le2\binom{r}{t}/\binom{s}{t}$ and $1\le\tau\le \binom{r}{t}/\binom{s}{t}$.
As $n\to\infty$,
\[
\begin{aligned}
 R_3(P_n^{\mathrm{lin}};r,s,t)&=\gamma n+o(n),&
 R_3(C_n^{\mathrm{lin}};r,s,t)&=\gamma n+o(n),\\
 R_3(P_n^{\mathrm{tight}};r,s,t)&=\tau n+o(n),&
 R_3(C_n^{\mathrm{tight}};r,s,t)&=\tau n+o(n)\quad(3\mid n).
\end{aligned}
\]
The linear coefficient satisfies $\gamma=2$ if and only if $r=s$.
\end{theorem}

We finally fix the length of a linear path and increase the number of
colors. For $h=\lfloor(n-1)/2\rfloor$ the largest number of edges in a
$P_n^{\mathrm{lin}}$-free triple system on $N$ vertices is
$(h+o(1))\binom N2$ as $N\to\infty$
\cite{pa:KMV,cy:JPR}. Applying this bound to each support and summing
gives $R_3(P_n^{\mathrm{lin}};r,s,t)\le(3h+o(1))\binom{r}{t}/\binom{s}{t}$.
Bohman and Zhu~\cite{pa:BZ} obtained a better coefficient for the ordinary
three-edge path using the directed-triangle theorem of Hladk\'y, Kr\'al'
and Norin~\cite{psix:HKN}. The next theorem extends their coefficient
to arbitrary fixed lengths, with a factor $h$ and with several colors
required in common.

\begin{theorem}\label{thm:growing-paths}
Fix $n\ge3$ and $s\ge t\ge2$, and let $h=\lfloor(n-1)/2\rfloor$.
As $r\to\infty$,
\[
 h\frac{\binom{r}{t}}{\binom{s}{t}}-O_{n,s,t}(r^{t-1})\ \le\ R_3(P_n^{\mathrm{lin}};r,s,t)\
 \le\left(\frac{2000h}{1307}+o(1)\right)\frac{\binom{r}{t}}{\binom{s}{t}}.
\]
In ordinary $k$-colorings the upper bound is
$(2000h/1307+o(1))k$ as $k\to\infty$.
\end{theorem}

The following example explains why the structure of the labels matters in
the proof of Theorem~\ref{thm:growing-paths}. Replacing each label $S$
by $\binom St$ turns $t$ common colors into one common element of
$\binom{[r]}t$. Allowing arbitrary subfamilies of $\binom{[r]}t$
of size $\binom{s}{t}$ can change the Ramsey number exponentially.

\begin{exampleintro}
For $m\ge1$, take $r=8m$, $s=4m$ and $t=m$. Then
$R_3(B_2;8m,4m,m)=4$. On four vertices, let $d_i$ count the triples
whose labels contain color $i$. Since $\sum_i d_i=2r$, the sum
of the six pairwise intersection sizes satisfies
$\sum_i\binom{d_i}{2}\ge((\sum_i d_i)^2/r-\sum_i d_i)/2=r$
by Cauchy--Schwarz. Some two triples therefore share at least $r/6>m$
colors and form $B_2$. The star itself needs four vertices.

For the relaxed labels, let
$k=\left\lfloor\binom{8m}{m}/\binom{4m}{m}\right\rfloor$. We have
\[
 k+1\le R_3\!\left(B_2;\binom{8m}{m},\binom{4m}{m},1\right)\le k+3.
\]
Indeed, choose $k$ disjoint families of $\binom{4m}{m}$ colors from
the $\binom{8m}{m}$ available colors. Index the families by
$\mathbb Z_k$ and label each triple $xyz$ on $\mathbb Z_k$ by the
family indexed by $x+y+z$. Distinct triples containing a fixed pair
have disjoint labels. Conversely, in an avoiding coloring on $N$
vertices the $N-2$ labels through any pair are disjoint, so
$(N-2)\binom{4m}{m}\le\binom{8m}{m}$. Each factor
$(8m-i)/(4m-i)$ is at least two for $0\le i<m$, so
$2^m\le\binom{8m}{m}/\binom{4m}{m}\le2^{8m}$. Hence the relaxed
Ramsey number grows exponentially in $m$.
\end{exampleintro}

\section{The proof of Theorem~\ref*{thm:stars}}\label{sec:stars}

\begin{theoremproof}
For a labeling $\phi$, a vertex $x$, and $T\in\binom{[r]}t$, let
$L_{x,T}$ be the \emph{link graph} on the other vertices whose edges are the pairs
$yz$ with $T\subseteq\phi(\{x,y,z\})$. A \emph{matching} is a set
of pairwise disjoint edges. Write $\nu(G)$ for its maximum size in
a graph $G$. A copy of $S_n$ with center $x$ and common colors $T$
is equivalent to a matching of size $n$ in $L_{x,T}$. Thus, if
$f(N)=\min_\phi\max_{x,T}\nu(L_{x,T})$, restricting a labeling to
fewer vertices shows that $f$ is nondecreasing. The first $N$
for which $f(N)\ge n$ is $R_3(S_n;r,s,t)$.

Let $\operatorname{ex}_3(v,S_n)$ be the largest number of edges
in an $S_n$-free $3$-graph on $v$ vertices. The theorem of
Chung and Frankl~\cite[Theorems~1.1 and~1.2]{ls:CF}, including
the two-edge case of Erd\H{o}s and S\'os recalled on page~111 of
that paper, gives $\operatorname{ex}_3(v,S_n)/v\to a_n$ for fixed
$n\ge2$. For a finite $S_n$-free $3$-graph with $v$ vertices and $e$ edges,
the union of $b$ disjoint copies is also $S_n$-free. Thus
$e/v\le\operatorname{ex}_3(bv,S_n)/(bv)\to a_n$, giving
$\operatorname{ex}_3(v,S_n)\le a_nv$ for every $v$. In an avoiding labeling on $v$ vertices, count pairs
consisting of a triple and a $t$-set contained in its label.
There are $\binom v3\binom st$ such pairs, and at most $a_nv$
triples for each $t$-set. Solving the resulting quadratic and using
the $2n+1$ vertices of the star gives
\[
 2n+1\le R_3(S_n;r,s,t)
 \le \left\lfloor\frac{3+\sqrt{1+24a_n\,\binom{r}{t}/\binom{s}{t}}}2\right\rfloor+1
 \le \sqrt{6a_n\,\binom{r}{t}/\binom{s}{t}}+3.
\]

First fix $r,s,t$. We estimate $f(N)$ by sampling the link graphs.
A \emph{fractional vertex cover} of a graph assigns a weight in
$[0,1]$ to each vertex so that the endpoints of every edge have
total weight at least one. Its weight is the sum of all vertex
weights and bounds the size of every matching. It is
\emph{half-integral} if all its entries lie in $\{0,1/2,1\}$.

For $2\le k\le v$, define
\(
\Delta(v,k)=\frac12\max_{\substack{1\le a\le v\\a\ {\rm odd}}}
 \max\left\{0,\ \frac1a-
        \frac{\binom{v-a}{k-1}}{\binom{v-1}{k-1}}\right\},
\)
where a binomial coefficient is zero if the selection is impossible.

\begin{claim}\label{cl:link-sample}
If $X$ is a uniformly chosen $k$-subset of a $v$-vertex graph $G$,
then $G[X]$ has a half-integral cover of random weight $Y$ such that
\[
 {\mathbb E}Y\le (k/v)\nu(G)+k\Delta(v,k),\qquad
 {\mathbb P}(Y-{\mathbb E}Y>u)\le \exp(-2u^2/k)\quad(u>0).
\]
\end{claim}
\begin{claimproof}
The Tutte--Berge formula~\cite{ls:TB} gives a vertex set
$U$ such that $\nu(G)=|U|+\sum_C\lfloor|C|/2\rfloor$, where $C$
runs over the components of $G-U$. Give the sampled vertices in $U$
weight one. In each remaining component give its sampled vertices
weight $1/2$ if at least two were sampled, and zero otherwise.
These weights cover every sampled edge. For $p=k/v$ and a component
$C$ of size $a$, the probability of sampling exactly one vertex of
$C$ is $pa\binom{v-a}{k-1}/\binom{v-1}{k-1}$. Hence
\[
 {\mathbb E}Y-p\nu(G)
 =\frac12\sum_C\left(p\mathbf 1_{\{|C|\ {\rm odd}\}}
              -{\mathbb P}(|C\cap X|=1)\right)
 \le k\Delta(v,k).
\]
For the inequality, an odd component of size $a$ contributes at
most $pa\Delta(v,k)$, an even component contributes at most zero,
and the component orders sum to at most $v$. Exchanging one selected
vertex for another changes $Y$ by at most one: each addition or
removal changes the weight in its component by an amount in $[0,1]$
in absolute value, and the two changes have opposite signs.
Expose the sample in order. Conditional on any prefix, couple
two possible next vertices by transposing them in the remainder
of the sample. The final subsets differ by at most one exchange,
so the conditional expectations lie in an interval of length one.
Applying Hoeffding's exponential moment bound~\cite{ls:Hoeffding}
at each of the $k$ steps gives the stated tail estimate.
\end{claimproof}

We next describe the finite linear programs. For $m$ vertex types
$1,\ldots,m$, assign a label to every multiset of three types,
allowing repeated entries. For each type $i$ and each $T\in\binom{[r]}t$,
choose weights $0\le w_{i,T}(j)\le1$ on the types, with
$w_{i,T}(j)+w_{i,T}(\ell)\ge1$ whenever the label on $\{i,j,\ell\}$
contains $T$. For a fixed labeling, minimize a variable $z$
subject to these constraints and $\sum_j w_{i,T}(j)\le mz$ for
every $i,T$. Let $\mu_m$ be the smallest optimum over all such
labelings on exactly $m$ types, and set $\mu=\inf_{m\ge1}\mu_m$.
For fixed $m$ there are only finitely many labelings. Each program
admits the all-one rows and has objective bounded below by zero, so its
rational coefficients give an attained rational optimum.

Replacing each type by $b$ vertices and assigning labels according
to types gives $f(bm)\le bm\mu_m$, by lifting the covers and omitting
the center itself. Consequently
$\liminf_{N\to\infty}f(N)/(N-1)\le\mu$. To bound $\mu$ from above
using a sample, we must also label triples with repeated types.

\begin{claim}\label{cl:link-extension}
Suppose that each link of a labeling on $m$ distinct vertices has a
half-integral cover of weight at most $W$. The labeling extends to all
multisets of three of these types, with a half-integral cover of weight
at most $W+2^{\binom{r}{t}}+1$ in every row.
\end{claim}
\begin{claimproof}
Fix a center $i$. For each other type $j$, let $Z(j)$
be the set of $T$ for which $w_{i,T}(j)=0$. Call $j$ bad if
every $s$-set contains a member of $Z(j)$. Two bad types cannot
have the same $Z(j)$, since their triple with $i$ would have two
zero endpoints in a link edge. There are at most $2^{\binom{r}{t}}$ bad types.
Raise all their entries in the rows at $i$ to one and give $i$
itself weight one in those rows. For $j\ne i$, label $\{i,j,j\}$
by an $s$-set all of whose $t$-subsets have positive modified
weight at $j$ in the rows at $i$. Half-integrality covers the
repeated pair at center $i$, and the diagonal entry one covers
the pair at center $j$. Label $\{i,i,i\}$ arbitrarily.
\end{claimproof}

Now take an optimal labeling on $v+1$ vertices and sample $k+1$
of them uniformly, with $2\le k\le v$. Conditional on the inclusion
of any center, the other $k$ sampled vertices are uniform among its
$v$ link vertices. By Claim~\ref{cl:link-sample}, with positive
probability every sampled link has a cover of weight at most
\[
 W=k\bigl(f(v+1)/v+\Delta(v,k)\bigr)+u_k,\qquad
 u_k=\sqrt{k\log(2\binom{r}{t}(k+1))/2}.
\]
Indeed, any fixed center is selected with probability $(k+1)/(v+1)$.
Conditional on its selection, each of its $\binom{r}{t}$ rows exceeds the bound
with probability at most $\exp(-2u_k^2/k)$. The union bound over all
$v+1$ centers and all rows is therefore $(k+1)\binom{r}{t}\exp(-2u_k^2/k)=1/2$.
Choose such a sample. Its link matching numbers are at most $W$,
and Claim~\ref{cl:link-extension} gives
\[
 f(k+1)\le W,\qquad
 \mu\le\mu_{k+1}\le\frac{W+2^{\binom{r}{t}}+1}{k+1}.
\]
For fixed $k$, $\Delta(v,k)\to0$ as $v\to\infty$. In fact the
binomial ratio tends to one for each fixed $a$, while the positive
part for $a>L$ is at most $1/L$. Let $v\to\infty$ along a sequence
attaining the lower limit of $f(v+1)/v$. The last two inequalities give
\[
 f(k+1)\le k\mu+u_k,\qquad
 0\le\mu_{k+1}-\mu\le\frac{u_k+2^{\binom{r}{t}}+1}{k+1}.
\]

For the other direction take
$k=\lfloor v^{2/3}\log v/16\rfloor$ for large $v$. We show that
$\Delta(v,k)=O(v^{-1/3})$. For $1\le a\le v^{1/3}$, writing the
binomial ratio as a product and using
$\log(1-x)\ge-x/(1-x)$ gives
\[
 \log\frac{\binom{v-a}{k-1}}{\binom{v-1}{k-1}}
 \ge-\frac{(k-1)(a-1)}{v-k-a+2}\ge-\frac{2k(a-1)}v .
\]
The concave function $\log a-2k(a-1)/v$ is nonnegative at both
ends of $[1,v^{1/3}]$. Thus the ratio is at least $1/a$ on
that interval. Larger $a$ contribute at most $v^{-1/3}$ to the
maximum. Also $u_k/k=O(v^{-1/3})$. Substituting in the sampled
upper bound for $\mu_{k+1}$ yields
$f(v+1)\ge v\mu-O_{r,s,t}(v^{2/3})$. Together with the bound
on $f(k+1)$, this gives $f(N)=\mu(N-1)+O_{r,s,t}(N^{2/3})$.
An optimal labeling on $N$ vertices avoids $S_{f(N)+1}$. Applying
the counting bound to this labeling and using $a_j\le j^2$ gives
$(N-1)(N-2)\binom{s}{t}\le6(f(N)+1)^2\binom{r}{t}$. Thus $\mu\ge1/\sqrt{6\,\binom{r}{t}/\binom{s}{t}}$ by
division by $N^2$ and passage to the limit. The link order gives
$\mu\le1/2$. For $M=R_3(S_n;r,s,t)$ we have $f(M-1)<n\le f(M)$
and $M=O_{r,s,t}(n)$ by the finite bound. Substitution of the estimates
for both $f(M-1)$ and $f(M)$ gives
$\mu M-O(M^{2/3})\le n\le\mu M+O(M^{2/3})$. Hence with $c=1/\mu$
we obtain $R_3(S_n;r,s,t)=cn+O_{r,s,t}(n^{2/3})$. Also
\[
 0\le c-\frac1{\mu_{k+1}}
   \le 6\frac{\binom{r}{t}}{\binom{s}{t}}\,\frac{u_k+2^{\binom{r}{t}}+1}{k+1}.
\]
The right-hand side tends to zero with $k$, giving intervals of arbitrarily small width.

In the ordinary two-color case take $n-1=3u+v$, $0\le v\le2$.
Partition $7u+2v$ vertices into sets $X,Y,Z$ of sizes
$4u+v,2u+v,u$. Color a triple blue if its class multiset is $XXY$,
$XXZ$, or $XYY$, and red otherwise. Use the cover weights in Table~\ref{tab:link-covers}.

\begin{table}[htbp]
\centering
\caption{Link-cover weights on $X,Y,Z$, in that order, for each choice of center.}
\label{tab:link-covers}
\begin{tabular}{@{}lcc@{}}
\toprule
Center class & Red link & Blue link\\
\midrule
$X$ & $(1/2,0,1)$ & $(0,1,1)$\\
$Y$ & $(0,1/2,1)$ & $(1/2,1/2,0)$\\
$Z$ & $(0,1,1)$ & $(1/2,0,0)$\\
\bottomrule
\end{tabular}
\end{table}
Deleting the center type from any allowed triple gives two entries
summing to at least one in the appropriate row. With the stated
class sizes every row has weight at most $3u+v=n-1$. This proves
$c(2,1,1)\ge7/3$. For the upper bound, the two-color matching
theorem of Cockayne and Lorimer~\cite{star:CL} forces a monochromatic
matching of size $n$ on $3n-1$ vertices. Apply it to the link of
one vertex in a red--blue labeling on $3n$ vertices to obtain
$c(2,1,1)\le3$.

For the tight star $B_n$, let $D$ be the largest length of a list of $s$-sets in $[r]$,
allowing repetitions, in which each $t$-set occurs at most $n-1$
times. Counting incidences gives $D\le(n-1)\binom{r}{t}/\binom{s}{t}$. In a labeling
avoiding $B_n$ on $N$ vertices, the labels of the $N-2$ triples
through any fixed pair form such a list. Conversely, write a maximum
list as $S_0,\ldots,S_{D-1}$. If $D\ge3$, label the triple
$\{x,y,z\}$ on $\mathbb Z/D\mathbb Z$ by $S_{x+y+z\bmod D}$.
At any fixed pair the distinct third vertices select distinct
positions, so this labeling avoids $B_n$. If $D<3$ the order of
$B_n$ suffices for the same lower bound. Hence
$D+1\le R_3(B_n;r,s,t)\le D+3$.

Each $t$-set is contained in $C=\binom{r-t}{s-t}$ possible
labels. Taking $\lfloor(n-1)/C\rfloor$ copies of every $s$-set
produces a list of length at least $(n-1)\binom{r}{t}/\binom{s}{t}-\binom rs$, since
$\binom{r}{s}\binom{s}{t}=\binom{r}{t}C$. This proves the bounded error with colors fixed.
For fixed $s,t$ and $r\to\infty$, let $K$ be the maximum number of
$s$-sets with no $t$-set contained in two of them. R\"odl's
packing theorem~\cite{stars:Rodl} gives $K=(1-o(1))\binom{r}{t}/\binom{s}{t}$ for $s>t$,
and $K=\binom{r}{t}$ if $s=t$. Repeating this family $n-1$ times gives
$D\ge(n-1)K$, proving $R_3(B_n;r,s,t)\sim(n-1)\binom{r}{t}/\binom{s}{t}$.

It remains to prove the lower bound for $S_n$ when $n,s,t$ are fixed
and $r\to\infty$. Fix a nonempty $S_n$-free $3$-graph $J$ without
isolated vertices, with $w$ vertices and $p$ edges. The extremal
limit above allows $p/w$ to approach $a_n$ from below.

\begin{claim}\label{cl:decomposition}
There is an arithmetic progression such that, for every sufficiently
large $v$ in it, $K_v^{(3)}$ decomposes into copies of $J$ that can
be partitioned into at most $(w/(6p)+o(1))v^2$ classes. Within each
class the copies are vertex-disjoint.
\end{claim}
\begin{claimproof}
For $i=0,1,2$, let $d_i$ be the greatest common divisor of the
numbers of edges of $J$ containing its individual $i$-sets, so
$d_0=p$. The design theorem of Glock, K\"uhn, Lo, and
Osthus~\cite{ls:GKLO} provides a $J$-decomposition of the complete
$3$-graph on $\ell$ vertices for all sufficiently large $\ell$
satisfying $d_i\mid\binom{\ell-i}{3-i}$ for $i=0,1,2$.
Choose $\ell\equiv2\bmod(6p d_1d_2)$ large enough. The three necessary
divisibilities follow by writing $\binom{\ell-i}{3-i}$ as the product
of $3-i$ consecutive integers ending at $\ell-2$, divided by $(3-i)!$.
Name the
copies in its decomposition $J_1,\ldots,J_b$, where $b=\binom\ell3/p$.
Write $b_i=\binom{\ell-i}{3-i}$. The same theorem, with the complete
$3$-graph on $\ell$ vertices as the block, gives a family of
$\ell$-subsets covering every triple exactly once for all sufficiently
large $v\equiv\ell\bmod (6\operatorname{lcm}(b_0,b_1,b_2))$.
Indeed, the numerator polynomials for $\binom{v-i}{3-i}$ and
$\binom{\ell-i}{3-i}$ have congruent values modulo $(3-i)!\,b_i$,
so these binomial coefficients agree modulo $b_i$.

Independently in each such block, choose a uniformly random bijection
from $[\ell]$ to the block and transport the fixed inner decomposition.
For each $j\in[b]$, the vertex sets of the transported copies of $J_j$
form a simple $w$-uniform hypergraph, since distinct outer blocks
meet in at most two vertices. Each vertex degree is a sum of independent
indicators of mean in total
$D_v=(w/\ell)\binom{v-1}2/\binom{\ell-1}2=\Theta(v^2)$.
Chernoff's bound and a union bound over $bv$ degrees show that
all are $D_v+O(v^{3/2})$ for some choice of the bijections. Any
pair belongs to at most $(v-2)/(\ell-2)=o(D_v)$ blocks. The
Pippenger--Spencer edge-coloring theorem~\cite{ls:PS} therefore
partitions each of these $w$-uniform hypergraphs into $(1+o(1))D_v$
matchings. Giving distinct classes to different $j$ yields
$b(1+o(1))D_v=(w/(6p)+o(1))v^2$ classes in all.
\end{claimproof}

Fix $0<\epsilon<1$. Given $q$ ordinary colors, take the largest $v$
in the progression of Claim~\ref{cl:decomposition} below
$(1-\epsilon)\sqrt{6pq/w}$. For large $q$ there are at most $q$
classes, so assign a different color to each class. Any
monochromatic $S_n$ would lie within one of its vertex-disjoint copies
of $J$, which is impossible. Taking $q=K$, replace each used color
by a distinct label from the family of $K$ labels above. Since any
$t$-set occurs in at most one of them, the resulting labeling also
avoids $S_n$. The progression has bounded gaps and $K\binom{s}{t}/\binom{r}{t}\to1$.
Letting $r\to\infty$ first, then $\epsilon\downarrow0$ and $p/w$
approach $a_n$, gives the required lower coefficient $\sqrt{6a_n}$.
\end{theoremproof}

\section{The proof of Theorem~\ref*{thm:paths-cycles}}\label{sec:paths-cycles}

\begin{theoremproof}
To determine the leading coefficients \(\gamma\) and \(\tau\), we use
two finite packing problems. The first applies to linear paths and cycles,
where consecutive edges share one vertex, and the second to tight paths
and cycles, where consecutive edges share two vertices.

Two edges of a \(3\)-graph are in the same \emph{component} if a chain of
edges with consecutive members sharing a vertex joins them. They are in
the same \emph{tight component} if the chain can be chosen with consecutive
members sharing a pair. The \emph{pair shadow} of a \(3\)-graph is the graph of pairs contained in its edges. We only count components containing edges.
Distinct components of the first kind have disjoint vertex sets. Distinct
tight components have disjoint pair shadows and may share vertices.

For a \(3\)-graph \(C\), define \(\nu_{211}(C)\) as follows. For each edge
\(abc\), allow a nonnegative weight with one vertex, say \(a\),
distinguished. One unit of this weight uses two units of capacity
at \(a\) and one each at \(b,c\). Each vertex has capacity one.
The maximum total weight is \(\nu_{211}(C)\). These capacities correspond
to the two triples \(a_1bc,a_2bc\), with \(a_1,a_2\) of the same type.
Fractional packings of this four-vertex configuration also appear in
work of Han, Sun and Wang~\cite{pa:HSW}. A \emph{fractional matching}
uses one unit of capacity at each of the three vertices per unit of
edge weight. Write \(\nu^*(C)\) for its maximum weight. Linear programming
duality expresses each maximum as a minimum of \(\sum_v w_v\) over
nonnegative vertex weights, called a \emph{cover}. For \(\nu_{211}(C)\),
require \(2w_a+w_b+w_c\ge1\) on each edge \(abc\) for every choice
of its distinguished vertex \(a\). For \(\nu^*(C)\), require
\(w_a+w_b+w_c\ge1\) on each edge. The weights may be truncated at
one in both cases. Adding the capacity
constraints gives \(\nu_{211}(C)\le |V(C)|/4\) and \(\nu^*(C)\le |V(C)|/3\).
In the first program some optimal solution has at most \(|V(C)|\)
positive variables, since a basic feasible solution has at most as many
positive variables as there are vertex constraints.

For a coloring \(\chi\) on \(m\ge3\) vertices, set
\[
 \begin{aligned}
 p_{\mathrm{L}}(\chi)&=\frac2m\max_{T,\ C\text{ component of }H_T}\nu_{211}(C),&
 g_{\mathrm{L}}(m)&=\min_\chi p_{\mathrm{L}}(\chi),\\
 p_{\mathrm{T}}(\chi)&=\frac3m\max_{T,\ C\text{ tight component of }H_T}\nu^*(C),&
 g_{\mathrm{T}}(m)&=\min_\chi p_{\mathrm{T}}(\chi).
 \end{aligned}
\]
The minima range over complete \((r,s)\)-colorings, with value zero
for a component maximum over an empty family. A partial coloring assigns labels to only some triples. We also use \(p_{\mathrm{L}},p_{\mathrm{T}}\) for partial colorings,
taking the supports of the labeled triples only. The weights are confined
to one component, but connecting edges may have weight zero.

If a \(t\)-common linear path has \(\ell\) edges, give each edge
\(x_{i-1}y_ix_i\) weight \(1/2\), distinguishing \(y_i\).
Each distinguished \(y_i\) uses one unit of capacity, and each \(x_i\)
lies in at most two edges of weight \(1/2\). Thus the total weight
\(\ell/2\) is feasible in one component, and \(\ell\le m p_{\mathrm{L}}(\chi)\).
In a tight path each vertex meets at most three edges. Giving each
edge weight \(1/3\) yields \(\ell\le m p_{\mathrm{T}}(\chi)\).

The incidence identity also gives lower bounds for both packing
problems. In a component for the linear program, summing the three
cover inequalities on \(abc\) gives \(w_a+w_b+w_c\ge3/4\).
The degree of any vertex \(v\) in \(C\) is at most \(|\partial C|\),
where \(\partial C\) is its pair shadow, since the map \(e\mapsto
e\setminus\{v\}\) sends incident edges to different pairs.
Thus \(e(C)\le(4/3)\nu_{211}(C)|\partial C|\) in the linear case and
\(e(C)\le\nu^*(C)|\partial C|\) in the tight case.
Component shadows in each support are disjoint in both cases.
For each \(T\), the sum of the shadow sizes is at most \(\binom m2\).
Bounding the packing in each component by the maximum in the definition
of \(p_X(\chi)\), summing over the \(\binom{r}{t}\) choices of \(T\), and dividing
the resulting bound for \(\binom{s}{t}\binom m3\) by \(\binom m2\) yields
\begin{equation}\label{pc:eq:crude-pack}
 \frac{(m-2)\binom{s}{t}}{2m\binom{r}{t}}\le g_{\mathrm{L}}(m)\le\frac12,\qquad
 \frac{(m-2)\binom{s}{t}}{m\binom{r}{t}}\le g_{\mathrm{T}}(m)\le1 .
\end{equation}

\begin{claim}\label{cl:sample-pack}
Let \(X\) stand for \(\mathrm{L}\) or \(\mathrm{T}\). For all sufficiently large \(m\)
depending only on \(r,t\), the following holds. In a partial coloring
\(\phi\) on \(v\ge m\) vertices, suppose that the fraction \(\rho\) of
missing triples satisfies \(\rho\binom m3\le1/4\).
Some \(m\)-set induces a complete coloring \(\chi\) with
\(p_X(\chi)\le p_X(\phi)+3m^{-1/4}\).
Consequently \(g_X(m)\) has a limit \(g_X\) as \(m\to\infty\).
\end{claim}
\begin{claimproof}
Put \(\delta=m^{-8}\). In each support call a component large if its
pair shadow has at least \(\delta v^2\) pairs. Disjointness of the
shadows leaves at most \(\binom{r}{t}/(2\delta)\) large components in all supports.
A component with \(p\) shadow pairs has at most \(p^{3/2}\) edges.
Indeed, its edges give triangles in the shadow graph. If \(D\) is
the adjacency matrix and \(\lambda_i\) its eigenvalues, the triangle
count is \(\operatorname{tr}(D^3)/6
\le(\sum_i\lambda_i^2)^{3/2}/6=(2p)^{3/2}/6\).
Summing over small components shows that at most
\(\binom{r}{t}\sqrt\delta\,v^3/2\) triples occur in any of them.

Choose a uniform \(m\)-set. The probability of seeing a triple from
a small component is at most \(2\binom{r}{t}/m\), by multiplying the preceding
bound by \(\binom m3/\binom v3\) and using \(v^3/\binom v3\le12\)
for \(v\ge6\). The probability of seeing any missing triple is at
most \(1/4\). For each large component, take an optimal cover in
\([0,1]\), extended by zero on the other vertices. Write \(a=2\)
in the linear case and \(a=3\) in the tight case. Its sampled weight
has mean at most \(mp_X(\phi)/a\). Hoeffding's inequality for sampling
without replacement \cite[Section~6]{ls:Hoeffding} bounds the probability
of exceeding this by \(m^{3/4}\) by \(\exp(-2\sqrt m)\).
The sum of the three failure bounds is at most
\(2\binom{r}{t}/m+1/4+(\binom{r}{t}m^8/2)\exp(-2\sqrt m)<1\) for large \(m\).

On a successful sample, each component belongs to a single original
large component of the same kind. Its restricted cover is feasible,
which proves the sampling bound. Apply this to a complete coloring
attaining \(g_X(v)\). Then \(g_X(m)-3m^{-1/4}\le g_X(v)\) for every
\(v\ge m\). Taking first \(\liminf_v\) and then \(\limsup_m\) gives
convergence.
\end{claimproof}

To realize a fixed linear cycle length, we use regular triples of
vertex sets. For disjoint sets \(U,V,W\), the \emph{density} of a \(3\)-graph across
them is the fraction of triples in \(U\times V\times W\) that are edges.
The triplet is \emph{\(\varepsilon\)-regular} if this density changes by
at most \(\varepsilon\) whenever each set is replaced by a subset of at
least an \(\varepsilon\) fraction of its size. Parts in a vertex partition
are called \emph{clusters}. A \(3\)-graph on their indices recording
selected cluster triples is a \emph{reduced \(3\)-graph}.

\begin{claim}\label{cl:linear-embed}
For \(0<d\le1\) and \(0<\eta<1/2\) there is \(\varepsilon>0\)
such that the following holds for every \(k\) and all sufficiently
large \(b=b(k,d,\eta)\). Let \(H\) have \(k\) disjoint clusters of
size \(b\), and let each edge of a reduced \(3\)-graph \(F\)
represent an \(\varepsilon\)-regular triplet of \(H\) of density at
least \(d\). For any component \(J\) of \(F\), the hypergraph \(H\)
contains a linear cycle of each integer length \(\ell\) with
\(\eta kb\le\ell\le2b\nu_{211}(J)-\eta kb\).
\end{claim}
\begin{claimproof}
Set \(\theta=\eta/32\) and choose
\(0<\varepsilon<d\theta^2/1000\). In any triplet under consideration,
three subsets of size at least \(\varepsilon b\) have density greater
than \(d/2\). A consequence we will use repeatedly is that, for two
fixed such subsets \(Y,Z\), fewer than \(\varepsilon b\) vertices
of the third cluster have fewer than \((d/2)|Y||Z|\) pairs of
neighbors in \(Y\times Z\). Otherwise the set of these vertices
would contradict the density condition.

Take a basic optimal packing on \(J\) and discard its weights below
\(\theta\). At most \(k\) weights were positive, so the discarded sum
is at most \(\theta k\). For a remaining weight \(\lambda\) on \(ijz\)
with \(i\) distinguished, set \(u=\lfloor(1-\theta)\lambda b\rfloor\).
Allocate a set \(A_0\) of \(2u\) vertices in cluster \(i\), a set
\(B_0\) of \(u\) vertices in cluster \(j\), and a set \(C_0\) of \(u\)
vertices in cluster \(z\). Make all allocations
disjoint. The capacity constraints use at most \((1-\theta)b\) vertices
in each cluster, so also choose a disjoint reserve \(R_i\) of
\(\lfloor\theta b\rfloor\) unused vertices in each cluster \(i\).
For sufficiently large \(b\), each \(u\ge\theta b/2\).
If no weight survives, then \(\nu_{211}(J)\le\theta k\) and the asserted
interval is empty. Otherwise write \(p\) for the number of allocations.

To connect the allocations, order them cyclically. The entrance of each will be in
its \(B_0\) cluster and its exit in its \(C_0\) cluster. From each
exit cluster to the next entrance cluster choose a shortest path in
the pair shadow of \(J\). If they agree, use a two-step walk through
a neighbor. Each step goes between two different clusters \(i,j\)
in a chosen edge \(ijz\) of \(J\). The total number \(q\) of steps
satisfies \(1\le q\le pk\le k^2\). A second choice replaces one step
\(i,j\) by \(i,z,j\), using the same reduced edge \(ijz\) for both
steps. We realize each step by three hyperedges below. The two choices
therefore use exactly \(3q\) and \(3q+3\) connecting edges.

Fix one of the route choices. An endpoint \(x\) in cluster \(i\) is
suitable for a step using \(ijz\) if it has at least
\((d/2)|R_j||R_z|\) pairs in its link in \(R_j\times R_z\).
Fewer than \(\varepsilon b\) vertices fail, by the density observation.
In each \(B_0,C_0\), exclude these unsuitable endpoints for the single
route step incident there. The remaining sets will be used for the
entire internal path, so both endpoints will meet their route requirements.

Within each allocation, fix an odd integer \(1\le L\le2(1-\theta)u\).
We construct a path with edges \(x_{j-1}a_jx_j\) for \(1\le j\le L\).
The \(a_j\) are distinct vertices of \(A_0\), and the \(x_j\) are
distinct eligible vertices alternating between \(B_0\) and \(C_0\),
starting in \(B_0\) and ending in \(C_0\).
Write \(A'\) for the unused part of \(A_0\). Start in \(B_0\) with
a vertex having at least \((d/2)|A'||C'|\) link pairs into \(A'\)
and the available set \(C'\) in \(C_0\). Such a vertex exists by
the density observation. Mark it used.

At any step let \(Y\) be the unused eligible vertices on the opposite
side, and \(Z\) those left on the current side. Maintain a last vertex
with at least \((d/4)|A'||Y|\) pairs in \(A'\times Y\). Up to the
desired length, \(|Y|,|Z|\ge\theta u-\varepsilon b-2\ge\theta^2 b/4\)
and \(|A'|\ge2\theta u\), for large \(b\). Fewer than \(\varepsilon b\)
vertices of \(Y\) have fewer than \((d/2)|A'||Z|\) pairs into
\(A'\times Z\). Since \((d/4)|Y|>\varepsilon b\), the link of the
last vertex contains a pair \(a'y\) with \(y\) outside this exceptional
set. Append the resulting edge. Removing \(a'\) costs at most \(|Z|\)
pairs at \(y\), so \(y\) satisfies the maintained bound for the next
step when \(b\) is large. Stop after the fixed odd number of
steps. The initial and terminal vertices are suitable for their
respective route steps.

To join the internal paths, choose a distinct reserve vertex at each
intermediate occurrence of a cluster on a route. This vertex must be
suitable for both incident steps. At most
\(2\varepsilon b\) candidates fail the two tests. The \(O(k^2)\)
fixed vertices can therefore be chosen distinct. A single step
with endpoints \(x\in V_i\), \(y\in V_j\) can now be realized by the
three edges
\[
 xwz_1,\qquad wvz_0,\qquad vyz_2,
 \qquad w\in R_j,\quad v\in R_i,\quad z_0,z_1,z_2\in R_z .
\]
Indeed, at least \(d|R_j|/4\) choices of \(w\) each have at least
\(d|R_z|/4\) choices for \(z_1\). Otherwise the total number of pairs
in the link of \(x\) would be less than \((d/2)|R_j||R_z|\).
A corresponding
set of \(v\)'s works for \(y\). Exclude all fixed endpoints and
all previously used connector interiors from these sets and from \(R_z\).
This costs \(O(k^2)\) vertices. All three remaining sets have size
at least \(\varepsilon b\), so they contain a middle edge \(wvz_0\).
Choose \(z_1,z_2\) from the indicated link neighborhoods, avoiding
all forbidden vertices and each other. Each step now meets the rest
of the construction only at its fixed endpoints. In particular,
the two steps of the detour through \(z\) use distinct reserve vertices
apart from their common endpoint. This remains true when a cluster
or a reduced edge is visited more than once.

It remains to choose the internal lengths so that the cycle has
exactly \(\ell\) edges. In allocation \(j\), let \(U_j\) be the largest odd integer at most
\(2(1-\theta)u_j\). Its internal length can be any of
\(1,3,\ldots,U_j\), for either choice of routes. We have
\(U_j\ge2(1-\theta)u_j-2\) and
\(\sum_j u_j\ge(1-\theta)b(\nu_{211}(J)-\theta k)-k\). Hence
\[
 \sum_{j=1}^p U_j
 \ \ge\ 2(1-\theta)^2b(\nu_{211}(J)-\theta k)-4k
 \ \ge\ 2b\nu_{211}(J)-3\theta kb-4k ,
\]
where we used \(\nu_{211}(J)\le k/4\). For large \(b\) this exceeds
\(2b\nu_{211}(J)-\eta kb\), while
\(p+3(q+1)\le k+3(k^2+1)<\eta kb\).

Now fix \(\ell\) in the asserted interval. Choose
\(c\in\{3q,3q+3\}\) with \(c\equiv\ell-p\pmod2\) and use the
corresponding routes. The bounds just proved give
\[
 0\le\frac{\ell-c-p}{2}\le\sum_{j=1}^p\frac{U_j-1}{2}.
\]
Distribute the integer \((\ell-c-p)/2\) among the allocations in
order, giving each the smaller of its capacity \((U_j-1)/2\) and the
number of units left. The upper bound ensures that all units are used.
If allocation \(j\) receives \(z_j\) units, take an internal path of length
\(1+2z_j\), then join the endpoints along the chosen routes.
Internal paths have disjoint allocations. All connector interiors
are fresh, and each fixed endpoint lies on exactly the two
edges consecutive there. Thus the cyclic sequence is a linear
cycle, with precisely \(c+\sum_j(1+2z_j)=\ell\) edges. Only the
lower bound on \(b\) depends on \(k\).
\end{claimproof}

We now pass from a coloring of triples to a coloring of cluster
triples. Call a length \emph{admissible} for linear cycles if it is at
least three, and for tight cycles if it is at least six and divisible
by three.

\begin{claim}\label{cl:reduce}
Fix \(X\in\{\mathrm{L},\mathrm{T}\}\), \(m\ge3\), and \(0<\eta,\rho<1/2\).
For every sufficiently large \(N\), any coloring of \(K_N^{(3)}\)
gives \(k\ge m\) clusters of equal size \(b\), covering at least
\((1-\eta)N\) vertices, and a partial \((r,s)\)-coloring \(\phi\)
of the cluster triples with at most \(\rho\binom k3\) missing.
Every admissible integer \(\ell\) in the interval
\[
 \eta kb\le\ell\le (p_X(\phi)-\eta)kb
\]
is the length of a \(t\)-common cycle of the indicated kind in the
original coloring.
\end{claim}
\begin{claimproof}
Write \(M=\binom rs\) for the number of possible labels and fix
\(d=1/(2M)\).

For linear cycles, apply the weak hypergraph regularity lemma \cite{ref:ChungWeak,ref:KNRS} to
the \(\binom{r}{t}\) supports simultaneously. In the form used here, for any
precision \(\varepsilon>0\) and any minimum number of clusters \(m\),
it gives, for large \(N\), \(m\le k\le K(\varepsilon,m,r,t)\) equal
clusters and at most \(\varepsilon N\) exceptional vertices. For each
support all but at most \(\varepsilon\binom k3\) cluster triples are
\(\varepsilon\)-regular. The simultaneous form follows from the usual
proof by summing the mean squared densities of the \(\binom{r}{t}\) edge indicator
functions at each partition refinement. Choose \(\varepsilon\) smaller
than \(\eta,\rho/\binom{r}{t}\) and the precision in Claim~\ref{cl:linear-embed}
for \(d,\eta\). This choice is independent of the eventual value of \(k\).

On each triplet regular for every support, assign an original \(s\)-label
occurring on at least a \(1/M\) fraction of its triples. Leave the other
triplets unlabeled. Each edge of a reduced support then represents a
regular triplet of density at least \(1/M\) in the corresponding original
support. Take a reduced component attaining \(p_{\mathrm{L}}(\phi)\) and apply
Claim~\ref{cl:linear-embed}. The cluster size threshold holds uniformly
for all \(k\le K(\varepsilon,m,r,t)\) by increasing \(N\).

For tight cycles, regularity must also test subsets of triples selected by pairs.
We describe the regular slice and cycle embedding lemmas of Allen et
al.~\cite[Lemmas~6 and~9]{cyextra:ABCM} in the form needed here. A \emph{regular
slice} on equal clusters is an auxiliary graph between the clusters.
All cluster pairs have the same density parameter \(d_2>0\).
Between any two such clusters, the fraction of present pairs in
\(U\times W\) differs from \(d_2\) by at most \(\varepsilon\) whenever
\(|U||W|>\varepsilon b^2\). On three clusters,
the three bipartite graphs form a \emph{triad}. The triples tested
are the triangles of this triad, and the \emph{relative density} of
a \(3\)-graph is the fraction of these triangles that are its edges.
At precision \(\varepsilon_3\), regularity on a triad means that a
single density parameter approximates, within \(\varepsilon_3\), the
density on each union of the triangle sets of a fixed number of
subgraphs of the triad, whenever that union has more than an
\(\varepsilon_3\) fraction of the triad's triangles.

Apply the regular slice lemma to the \(M\) edge-disjoint systems
consisting of triples with each exact label. It gives a common slice
with \(m\le k\le K\) in which all but \(\varepsilon_3\binom k3\) triads
are regular for each exact label. The parameters can be chosen in the
following order. Fix the embedding loss \(\eta\), the density threshold
\(d\) and the missing fraction \(\rho\). Choose \(\varepsilon_3\) so
that \(M\varepsilon_3<\min\{\rho,1/(2M)\}\) and the cycle embedding lemma applies
with top-level error \(M\varepsilon_3\). Supply the slice lemma with
a lower-level precision and a number of test subgraphs as functions of
an upper bound \(u\) on the cluster count. For each \(u\) choose them,
respectively, sufficiently small and sufficiently large for every
allowed \(d_2\ge1/u\) with \(1/d_2\) an integer. The slice lemma then
returns \(K\) and a threshold on \(N\). Its divisibility assumption
is met by first discarding fewer than \(K!\) vertices. The equal
clusters therefore cover \(kb=N-O(K!)\) vertices. The lower-level
precision ensures that each triad contains triangles by the counting
estimate \cite[Fact~3]{cyextra:ABCM}. Once \(K\) is fixed there are
only finitely many possibilities for \(k\) and for \(1/d_2\), both
at most \(K\). Increase \(N\) so the embedding lemma's order
threshold holds for all of them.

On each triad regular for all \(M\) classes choose a label of relative
density at least \(1/M\). Such a label exists because the exact classes
partition its triangles. These labels give \(\phi\), with at most
\(M\varepsilon_3\binom k3\) missing triples. Every original \(H_T\)
is the disjoint union of the exact classes whose labels contain \(T\).
Their regularity discrepancies on any common test set add. Hence all
good triads are regular for \(H_T\) at the enlarged precision
\(M\varepsilon_3\), with the sum of the density parameters truncated
at one if necessary. The omitted fraction is at most the same
\(M\varepsilon_3\). Thus the common slice is a regular slice for
each union \(H_T\) with these parameters. On a triad labeled by a superset of \(T\), the
actual relative density of \(H_T\) is at least \(1/M\). Testing the full
triad shows that its regularity density parameter is at least
\(1/M-M\varepsilon_3>d\).

If \(p_{\mathrm T}(\phi)>0\), take a color set \(T\) and a tight
component \(C\) of its reduced support attaining the maximum. Every
edge of \(C\), including those of zero matching weight, is a regular
triad of \(H_T\) of relative density at least \(d\). The chains
joining the positive-weight edges of a maximum fractional matching
in \(C\) therefore also join them in the \(d\)-reduced graph of
\(H_T\). Extend this matching by zero outside \(C\). Its weight is
\(\omega=kp_{\mathrm T}(\phi)/3\). Apply the cycle embedding lemma
to this single support on \(kb\) vertices with \(k\) clusters.
It gives a tight cycle in \(H_T\) for each length \(\ell\) divisible
by three up to
\[
 (1-\eta)\,3\omega\,\frac{kb}{k}
  =(1-\eta)p_{\mathrm T}(\phi)kb
  \ge (p_{\mathrm T}(\phi)-\eta)kb ,
\]
since \(p_{\mathrm T}(\phi)\le1\). The edges of each such cycle
share \(T\). For zero packing value the asserted interval is empty.
\end{claimproof}

Fix \(X=\mathrm{L}\) or \(\mathrm{T}\) and let \(g=g_X>0\) be the limit from
Claim~\ref{cl:sample-pack} and \eqref{pc:eq:crude-pack}. On \(N\)
vertices a coloring attaining \(g_X(N)\) has no \(t\)-common path
of the indicated kind longer than \(Ng_X(N)\). For any \(0<\xi<1/g\),
take \(N=\lfloor(1/g-\xi)n\rfloor\). Since \(Ng_X(N)/n\to1-g\xi<1\),
the longest such path has fewer than \(n-2\) edges for all sufficiently
large \(n\). Deleting one
edge of a linear \(n\)-cycle yields a path of length \(n-1\);
deleting the two wrap-around edges of a tight \(n\)-cycle yields a
path of length \(n-2\). Hence the coloring avoids both the path
and the cycle of length \(n\). Letting \(\xi\downarrow0\) shows that
both lower limits of the Ramsey ratios are at least \(1/g\).

For the upper bound take \(N=\lceil(1/g+\xi)n\rceil\) with \(\xi>0\)
and put \(a=(1/g+\xi)^{-1}<g\). Choose \(\eta>0\) so small that
\(\eta<a/2\) and \((g-3\eta)(1-\eta)>a\). Fix a sufficiently large
\(m\) with \(g_X(m)-3m^{-1/4}>g-\eta\), then choose
\(\rho\binom m3\le1/4\) with \(\rho>0\). Apply Claim~\ref{cl:reduce}
and then Claim~\ref{cl:sample-pack} to the resulting reduced coloring.
Some complete \(m\)-set has packing value at most
\(p_X(\phi)+3m^{-1/4}\), so \(p_X(\phi)\ge g_X(m)-3m^{-1/4}>g-\eta\).
Since \(kb\le N\) and \(\eta<a/2\), the lower endpoint \(\eta kb\)
is below \(n\) for large \(n\). The upper endpoint is greater than
\((g-2\eta)(1-\eta)N\); its coefficient exceeds \(a\), whereas
\((n+4)/N\to a\). Thus \([n,n+4]\) lies in the embedding interval.
Choosing \(\ell=n\) gives a linear cycle for all large \(n\), and a
tight cycle for all large \(n\) with \(3\mid n\). To get a linear
path of length \(n\), delete one edge from a cycle of length \(n+1\).
For a tight path, choose \(\ell\) to be the smallest multiple of
three at least \(n+2\), so \(\ell\le n+4\). The first \(n+2\)
distinct vertices in its cyclic order give a tight path of \(n\) edges.
Letting \(\xi\downarrow0\) identifies \(\gamma=1/g_{\mathrm{L}}\) and
\(\tau=1/g_{\mathrm{T}}\). Equation~\eqref{pc:eq:crude-pack} gives the claimed
upper endpoints.

If \(r=s\), every label is \([r]\) and
\(R_3(P_n^{\mathrm{lin}};r,s,t)=2n+1\), so \(\gamma=2\).
Now suppose \(r>s\). For each \(i\in[r]\) fix an \(s\)-element label
omitting \(i\). Partition the vertices into ordered classes
\(V_1,\ldots,V_r\) of sizes \(m,6m,\ldots,6^{r-1}m\), where \(m\)
is a positive integer. Give each triple the label of the first class
it meets. Fix a \(t\)-element color set \(T\) and choose \(i\in T\).
Every edge of \(H_T\) meeting \(V_i\) contains a vertex of an earlier
class. In a linear path each earlier vertex is in at most two edges,
and each such edge contains at most two vertices of \(V_i\).
Thus at most \(4\sum_{j<i}|V_j|\) vertices of \(V_i\) lie on the
path. Since \(|V_i|-4\sum_{j<i}|V_j|\ge m\), every \(t\)-common
linear path misses at least \(m\) vertices. Write
\(v=1+6+\cdots+6^{r-1}\) and take \(m=\lfloor 2n/(v-1)\rfloor\).
Such a path uses at most \((v-1)m\le2n\) vertices, whereas
\(P_n^{\mathrm{lin}}\) has \(2n+1\). Hence
\(R_3(P_n^{\mathrm{lin}};r,s,t)>vm\), giving
\(\gamma\ge 2v/(v-1)>2\).
\end{theoremproof}

\section{The proof of Theorem~\ref*{thm:growing-paths}}\label{sec:growing-paths}

\begin{theoremproof}
Fix \(s\ge t\ge2\) and \(n\ge3\), and let \(r\to\infty\).
With \(h=\lfloor(n-1)/2\rfloor\), we have \(n=2h+1\) or \(n=2h+2\).
For the lower bound, let \(q\) be the maximum number of \(s\)-element labels whose pairwise intersections have size less than \(t\). Then \(\binom{r}{t}/\binom{s}{t}-O_{s,t}(r^{t-1})\le q\le \binom{r}{t}/\binom{s}{t}\).
The upper bound counts the \(\binom{s}{t}\) different \(t\)-sets in each label.
If \(s=t\) take all labels. Otherwise we use a \emph{\(t\)-design
of index one}, a family of \(s\)-sets in which every \(t\)-set occurs
exactly once. By Keevash's existence theorem \cite{cy:Keevash}, such a
design exists on every sufficiently large \(r'\) satisfying
\(\binom{s-i}{t-i}\mid\binom{r'-i}{t-i}\) for \(0\le i<t\).
To find one with \(r'=r-O_{s,t}(1)\), take the largest \(r'\le r\)
congruent to \(s\) modulo \(t!\prod_{i=0}^{t-1}\binom{s-i}{t-i}\).
For \(j=t-i\), the polynomial \(j!\binom{x-i}{j}\) has integer
coefficients, and the modulus is divisible by \(j!\binom{s-i}{j}\).
This gives each required divisibility. The design supplies
\(\binom{r'}t/\binom{s}{t}=\binom{r}{t}/\binom{s}{t}-O_{s,t}(r^{t-1})\) such labels on \([r]\).
Take \(q\) such labels and \(q\) ordered vertex classes. Give the
first \(q-1\) classes \(h\) vertices each and the last class \(2n\)
vertices. Assign to each triple the label of the first class it meets.
Since two distinct labels have fewer than \(t\) common colors, all
edges of a \(t\)-common path must have the same label. For any of the
first \(q-1\) labels its edges would all meet the corresponding class.
Each of its \(h\) vertices can meet at most two path edges, and
\(2h<n\). For the last label the path must lie on \(2n\) vertices,
one fewer than \(P_n^{\mathrm{lin}}\) requires. This proves
\(R_3(P_n^{\mathrm{lin}};r,s,t)\ge 2n+1+h(q-1)
\ge h\binom{r}{t}/\binom{s}{t}-O_{n,s,t}(r^{t-1})\).

For the upper bound, we first describe a single path-free support
using pairs of high \emph{codegree}, where \(d_H(uv)\) counts the
vertices \(z\) with \(uvz\in E(H)\). The constants in the next three
claims are independent of the choice of support.

\begin{claim}\label{cl:gap}
There is a constant \(C_n\) such that any \(N\)-vertex
\(P_n^{\mathrm{lin}}\)-free system has at most \(C_nN\) pairs of codegree
at least \(h+1\). For fixed \(0<\delta<1\) and sufficiently large
\(N\), deleting at most \(C_n\delta N^2+O_n(N)\) edges gives a
system \(F\) and a graph \(J\) of at most \(C_nN\) pairs such that
each edge of \(F\) contains a pair of \(J\). The pairs of \(J\)
have \(F\)-codegree at least \(\delta N\), and every other pair has
\(F\)-codegree at most \(h\).
\end{claim}
\begin{claimproof}
For every pair \(uv\) of codegree at least \(h+1\), choose \(h+1\)
third vertices. Include each vertex in a set \(X\) independently with
probability \(1/2\). Keep \(uv\) when both ends lie in \(X\) and
all its chosen third vertices lie outside \(X\). Put \(L=2^{2h+1}\).
The kept graph contains no ordinary path of \(L\) edges.
Indeed, the path expansion lemma of Kostochka, Mubayi and
Verstra\"ete \cite[Lemma~4.1]{pa:KMV} states that such a path of
pairs, each with \(h+1\) third vertices outside the whole path,
gives a linear path of \(n\) triples. By the graph path bound of
Erd\H{o}s and Gallai \cite{ls:EG}, at most \((L-1)|X|/2\) pairs were
kept. If \(z\) is the original number of eligible pairs, the expected
number kept is \(2^{-h-3}z\). Since \(\mathbb E|X|=N/2\), this yields
\(z\le 2^{h+1}(L-1)N\), as claimed.

Repeatedly delete all triples on a pair of current codegree in
\([h+1,\delta N)\). Each such pair had original codegree at least
\(h+1\), and it is chosen only once. The cost is at most
\(C_n\delta N^2\), where \(C_n=2^{h+1}(L-1)\). Let \(J\) consist of the pairs whose codegree
is now at least \(\delta N\). Triples with no pair of \(J\) form a
system of maximum pair codegree at most \(h\). If a subgraph of
this system had minimum vertex degree greater than \(2(n-1)h\),
one could grow a linear path to length \(n\). At an end vertex
of a path with \(j<n\) edges, at most \(2jh\) incident triples contain
any other previously used vertex. Thus one can repeatedly remove
a vertex of degree at most \(2(n-1)h\), which bounds the number of
these triples by \(2(n-1)hN\). Delete them. No \(J\)-pair loses
an extension in this last step.
\end{claimproof}

\begin{claim}\label{cl:centers}
For fixed \(n,\delta\) in Claim~\ref{cl:gap} and all sufficiently
large \(N\), there is a set \(B\) of \(O_{n,\delta}(1)\) centers.
All vertices outside \(B\) have bounded \(J\)-degree, and deleting
\(O_{n,\delta}(1)\) pairs from \(J-B\) leaves a matching and
isolated vertices. For \(b\in B\), let \(A_b\) be its neighbors
in \(J\) outside \(B\). These sets satisfy \(|A_b|\ge\delta N-O_{n,\delta}(1)\).
The graph \(I\) joining centers \(b,c\) when \(|A_b\cap A_c|\ge10n\)
has no path on \(h+1\) vertices.
\end{claim}
\begin{claimproof}
We first bound the number of vertex-disjoint ordinary two-edge paths
in \(J\). Take an integer \(u>\max(h,1/\delta)\) and let \(R\) be large
enough that every \(16\)-coloring of the pairs of \([R]\) has a
monochromatic \(u\)-set. Suppose \(R\) disjoint two-edge paths exist,
and set aside their \(3R\) vertices. Denote their edges by
\(e_i^1,e_i^2\). Let \(L_i^a\) be the set of third vertices extending
\(e_i^a\) in \(F\) outside the vertices set aside, so
\(|L_i^a|\ge\delta N-3R\). For \(i<j\) record the four bits specifying
whether \(|L_i^a\cap L_j^b|\ge10n\) for \(a,b\in\{1,2\}\).
On \(u\) indices in increasing order this \(2\)-by-\(2\) matrix is constant.

If an off-diagonal entry \((a,b)\) is one, traverse each of the
first \(h+1\) selected paths from its edge of type \(b\) to its
edge of type \(a\). If both diagonal entries are one, instead
alternate the two traversal orders. In either case consecutive
paths can be joined by choosing a common third vertex for the
exit edge of one and the entrance edge of the next. Choose all
\(h\) such vertices distinct and choose fresh third vertices for
the two terminal edges. This gives \(2h+2\) triples in a linear
path. Within a two-edge path the consecutive triples share its
middle vertex, and between paths they share their chosen third
vertex. It contains the forbidden length \(n\).
In the remaining case some diagonal entry \((a,a)\) is zero.
The union of the \(u\) sets \(L_i^a\) has size at least
\(u(\delta N-3R)-\binom u2(10n-1)>N\) for large \(N\),
since \(u\delta>1\). This is impossible.

The vertices of a maximal collection of disjoint two-edge paths
therefore form a bounded set \(S\) such that \(J-S\) has maximum
degree at most one. A vertex of \(J\)-degree \(D\) satisfies
\(\delta ND\le D(D-1)+h(N-1-D)\). To see this, sum the codegrees of its \(D\)
incident \(J\)-pairs, counting twice triples whose other two
vertices are both neighbors, and bound contributions with a
nonneighbor by the codegrees of the corresponding non-\(J\) pairs.
If \(D\le\delta N/2\), then \(D\le2h/\delta\). Otherwise
\(D\ge\delta N-2h/\delta\). Take \(B\) to be the vertices of degree
exceeding a constant larger than \(|S|+1\) and \(2h/\delta\).
They all belong to \(S\), which proves the assertion on their
neighborhood sizes. Outside \(B\), only a bounded number of pairs
are incident to \(S-B\), and the other pairs form a matching.

If \(b_1,\ldots,b_{h+1}\) formed a path in \(I\), insert a distinct
common neighbor between every two successive centers. The resulting
graph path in \(J\) has \(2h\) edges. Extend it by a new neighbor
at one end when \(n=2h+1\), and at both ends when \(n=2h+2\).
All pairs of the resulting path have at least \(\delta N\) third
vertices in \(F\). Choose one per pair, distinct and outside the
whole graph path. The resulting triples form the forbidden linear path.
\end{claimproof}

A \emph{digraph} here has no loops and at most one arc in each
direction between distinct vertices. A pair carrying both arcs is
\emph{bidirectional}. The \emph{outdegree} counts arcs leaving a
vertex, and the \emph{indegree} counts arcs entering it.

\begin{claim}\label{cl:representation}
For every \(\rho>0\), if \(N\) is sufficiently large, a
\(P_n^{\mathrm{lin}}\)-free system \(H\) has a subgraph \(H^*\) missing
at most \(\rho N^2\) of its edges and a digraph \(M\) of maximum
outdegree at most \(h\) with the following properties. Its arcs
point from outside a bounded set \(B\) into \(B\), or occur in
both directions on a matching outside \(B\). Each edge of \(H^*\)
contains a matching pair or has the form \(bxy\) with exactly one
center \(b\in B\) and both \(x\to b,y\to b\) in \(M\).
A pair with no arc in either direction has \(H^*\)-codegree at most
\(h\). The bound on \(|B|\) depends only on \(n,\rho\).
\end{claim}
\begin{claimproof}
Choose \(0<\delta<\rho/(2C_n)\) and use
Claims~\ref{cl:gap} and~\ref{cl:centers}. Each noncenter is in
\(O_{n,\delta}(N)\) edges of \(F\). An incident edge contains one
of its boundedly many incident \(J\)-pairs or one of the \(O(N)\)
other \(J\)-pairs. Delete all edges meeting \(A_b\cap A_c\) for two
centers nonadjacent in \(I\). The union of these intersections
has bounded size and consists of noncenters, so this costs \(O(N)\)
edges. Discard all corresponding incident \(J\)-pairs. The centers
adjacent to any remaining noncenter are now pairwise adjacent in \(I\),
so they form a clique
of size at most \(h\) by Claim~\ref{cl:centers}. Delete also every
edge with at least two centers and every edge containing one of the bounded
number of exceptional pairs of \(J-B\). The former deletions cost at
most \(\binom{|B|}{2}N\) edges and each exceptional pair costs at most
\(N-2\). Discard these pairs. The surviving noncenter pairs of \(J\)
form a matching. Every surviving triple still contains an undiscarded
\(J\)-pair, since all triples on each discarded pair were removed.

An endpoint \(x\) of a matching pair might still have \(h\) adjacent
centers. Let \(C=\{b_1,\ldots,b_h\}\) be these centers. They form a
clique in \(I\). For the next constructions work in the original \(F\)
of Claim~\ref{cl:gap}, where the common-neighbor sets still have size
at least \(10n\) and all \(J\)-pairs have \(\delta N\) extensions.
For \(n=2h+1\), take a triple on the matching pair with a third
vertex outside \(B\). Starting at \(x\), take the \(J\)-path
\[
 x,b_1,u_1,b_2,u_2,\ldots,u_{h-1},b_h,y ,
\]
where the \(u_i\) are fresh common neighbors and \(y\) is a fresh
neighbor of \(b_h\). For \(h=1\) this is \(x,b_1,y\).
Choose all new vertices outside the first triple and expand the path's
\(2h\) pairs with distinct fresh third vertices. Together with the
first triple this is a forbidden path.

For \(n=2h+2\), two distinct matching pairs with endpoints \(x,y\)
adjacent to the same \(C\) would give the same path, now with \(y\)
as its fixed last vertex. First choose a triple on each matching
pair, outside \(B\) and avoiding the other pair and the other chosen
third vertex. There are at least \(\delta N\) possible extensions per pair
and only a bounded number of forbidden vertices. Choose all intermediate
neighbors and expansion vertices outside these triples. The \(2h\)-edge
path then meets the two triples only at its respective endpoints
\(x,y\). Appending the two triples gives a forbidden path of \(2h+2\)
edges. Thus each \(h\)-set of centers accounts for at most one matching
pair with an endpoint adjacent to all of them. Delete all edges on
these at most \(\binom{|B|}{h}\) pairs, at a cost of \(O(N)\).

Put both arcs on each remaining matching pair and direct each
remaining center--noncenter \(J\)-pair toward its center. The maximum
outdegree is at most \(h\), since a matching endpoint has at most
\(h-1\) adjacent centers. Finally discard any triple containing no
matching pair whose unique center fails to receive an arc from one
of its two other vertices. Such a triple uses a center--nonneighbor
pair of original codegree at most \(h\) in \(F\), since all triples
on previously discarded pairs were deleted. This costs at most
\(|B|hN\) edges. A pair unrepresented by an arc was either outside
\(J\) or lost all its edges. The total cost is
\(C_n\delta N^2+O_{n,\delta}(N)\le\rho N^2\) for large \(N\).
\end{claimproof}

We next compare the directions assigned by different supports.
Fix \(\lambda>1\), take \(N=\lceil\lambda h\binom{r}{t}/\binom{s}{t}\rceil\), and suppose a
coloring avoids the required path. For each \(T\in\binom{[r]}t\)
apply Claim~\ref{cl:representation} with the same \(\rho>0\).
Write \(H_T^*,M_T,B_T\) for its subgraph, digraph and centers.
A host triple \(e\) is \emph{clean} if \(e\in H_T^*\) for every
\(T\in\binom{\chi(e)}t\). At most \(\rho \binom{r}{t}N^2=O(\rho N^3)\)
triples fail to be clean, since \(\binom{r}{t}/N=O_{n,s,t,\lambda}(1)\).

For an \(s\)-label \(S\), call \(\binom St\) its \emph{complete
\(t\)-shadow}. Define \(M\) to contain \(x\to y\) if there exists
a label \(S\) such that \(x\to y\) belongs to every \(M_T\) with
\(T\in\binom St\). Every outgoing arc in \(M\) uses at least \(\binom{s}{t}\)
outgoing arcs across the \(\binom{r}{t}\) individual digraphs, counted with
their support indices. Consequently
\begin{equation}\label{pc:eq:outgoing}
 d_M^+(x)\le h\frac{\binom{r}{t}}{\binom{s}{t}} \qquad\text{for every }x .
\end{equation}
Here \(d_M^+(x)\) denotes outdegree. The next count uses the fact
that the \(t\)-sets in a label form its complete \(t\)-shadow.

\begin{claim}\label{cl:shadow}
Let \(\mathcal D\subseteq\binom{[r]}t\) and let a multiset of
\(s\)-labels be given. If each \(t\)-set outside \(\mathcal D\)
occurs in at most \(h\) labels, counting repetitions, the number
of labels whose complete \(t\)-shadow is not contained in
\(\mathcal D\) is at most \(h\binom{r}{t}/\binom{s}{t}+ht|\mathcal D|\binom{r-1}{t-1}\).
\end{claim}
\begin{claimproof}
Labels whose shadows are disjoint from \(\mathcal D\) each use
\(\binom{s}{t}\) incidences outside it, so there are at most \(h\binom{r}{t}/\binom{s}{t}\).
For a label \(S\) whose shadow meets both sides, set
\(Z=\bigcup_{T\in\mathcal D}T\). There is a \(t\)-set in \(S\)
outside \(\mathcal D\) meeting \(Z\). If \(S\subseteq Z\), choose any
member of its shadow outside \(\mathcal D\). Otherwise, choose one
color in \(S\cap Z\) and one in \(S\setminus Z\), and extend them
inside \(S\) to a \(t\)-set. This is where \(t\ge2\) is used.
There are at most \(|Z|\binom{r-1}{t-1}\) possible witnesses, with
at most \(h\) labels per witness. Now use \(|Z|\le t|\mathcal D|\).
\end{claimproof}

A pair is \emph{uncovered} if neither direction occurs in \(M\).
An \emph{in-star} on three vertices consists of arcs \(x\to z,y\to z\).
In the next estimates the constants multiplying \(\rho\) depend only
on \(n,s,t,\lambda\). Terms denoted by \(o(1)\) tend to zero as
\(r\to\infty\) at fixed \(\rho,n,s,t,\lambda\).

\begin{claim}\label{cl:agreement}
There are at most \((O(\rho)+o(1))N^2\) uncovered pairs.
All but \((O(\rho)+o(1))N^3\) host triples contain an in-star of
\(M\) or a bidirectional pair of \(M\).
\end{claim}
\begin{claimproof}
For a pair \(xy\), let \(\mathcal D_{xy}\) consist of the \(T\) for
which \(M_T\) contains some direction on \(xy\). Counting arcs gives
\(\sum_{xy}|\mathcal D_{xy}|\le h\binom{r}{t}N\).
Thus all but \(O(N^2/\sqrt r)\) pairs have
\(|\mathcal D_{xy}|\le\sqrt r\). Fix a pair satisfying this bound. On the
clean triples through \(xy\) each \(T\notin\mathcal D_{xy}\) occurs
at most \(h\) times by the last assertion of
Claim~\ref{cl:representation}. Claim~\ref{cl:shadow} bounds the
number of labels whose shadows are not contained in \(\mathcal D_{xy}\)
by \(h\binom{r}{t}/\binom{s}{t}+O(N/\sqrt r)\), using
\(\binom{r-1}{t-1}=t\binom{r}{t}/r\). If \(\mathcal D_{xy}\) contains no
complete \(t\)-shadow, at least \(N-2-h\binom{r}{t}/\binom{s}{t}-O(N/\sqrt r)\) triples
through \(xy\) fail to be clean. This is a positive constant
fraction of \(N\), since \(\lambda>1\). For all sufficiently large
\(r\) this fraction is at least \((1-1/\lambda)/2\).
Counting incidences with the at most \(\rho \binom{r}{t}N^2\) unclean triples,
each incident with three pairs, bounds the number of these pairs
by \(6\rho \binom{r}{t}N/(1-1/\lambda)=O(\rho N^2)\). In particular this
implicit constant is independent of the deletion parameter \(\rho\).

Suppose next that \(\binom St\subseteq\mathcal D_{xy}\). For each
\(T\) in this shadow, the pair carries an arc only toward \(x\),
an arc only toward \(y\), or both arcs on a matching pair in \(M_T\).
If these types agree throughout, \(xy\) is covered
in \(M\). Otherwise two \(t\)-sets differing in just one element
have different types, since one can move between any two members
of \(\binom St\) by single-element replacements. There are
\(\binom{r}{t}t(r-t)/2\) unordered adjacent pairs of \(t\)-sets in the whole
palette. For a fixed such pair \(T,T'\), opposite one-way arcs
force one endpoint into \(B_T\) and the other into \(B_{T'}\), in
some order. This gives \(O_{n,\rho}(1)\) vertex pairs. For a one-way
arc and a matching pair the center of the arc has at most one
partner in the other matching, giving the same bound. Since
\(N\) has order \(r^t\), these exceptions total
\(O_{n,s,t,\lambda,\rho}(\binom{r}{t}r)=o(N^2)\) for \(t\ge2\).
This proves the assertion on pairs.

For the assertion on triples fix a clean \(e\). In each support
indexed by \(T\in\binom{\chi(e)}t\), record its matching pair if it
has one, and otherwise its center. There is at most one matching
pair in a triple and the other case gives a unique center. A
constant record gives, respectively, a bidirectional pair or two
arcs toward the same vertex in \(M\). With different records there
are again two adjacent \(t\)-sets with different records. For a
fixed adjacent \(T,T'\), two different centers give \(O_{n,\rho}(N)\)
triples, since both center sets are bounded. For a center and a
matching pair, either the center determines its unique partner or
the third vertex is the center and there are at most \(N/2\)
matching pairs. Two different matching pairs in one triple must
intersect, and their common vertex determines both partners. These
cases each give \(O_{n,\rho}(N)\) triples. Summing over the adjacent
\(t\)-sets gives \(O_{n,s,t,\lambda,\rho}(\binom{r}{t}rN)=o(N^3)\), since
\(t\ge2\). Add the \(O(\rho N^3)\) triples that failed to be clean.
\end{claimproof}

We finish the upper bound by counting triangles among the arcs of \(M\)
whose reverses are absent. These arcs form an \emph{oriented graph},
a digraph with at most one of the two opposite arcs on each pair. A \emph{directed triangle} has arcs \(x\to y\),
\(y\to z\), \(z\to x\). Set \(\alpha=0.3465\).
Hladk\'y, Kr\'al' and Norin \cite{psix:HKN} proved that minimum
outdegree at least \(\alpha\) times the order forces a directed
triangle in an oriented graph. Reversing every arc gives the
same statement for minimum indegree.

It follows that for every \(g>0\) there is
\(c_g>0\) such that every sufficiently large oriented graph of order
\(v\) and minimum indegree at least \((\alpha+g)v\) has at least
\(c_gv^3\) directed triangles. Indeed, take a uniform subset of
\(u\) vertices, with \(u\) a sufficiently large constant depending
on \(g\). Conditional on any one of its vertices, its sampled
indegree has mean at least \((\alpha+g)(u-1)\).
For large \(u\), this mean exceeds \(\alpha u\) by at least \(gu/2\).
Hoeffding's inequality for the \(u-1\) draws without replacement
\cite[Section~6]{ls:Hoeffding} bounds each conditional failure
probability by \(\exp(-g^2u/2)\). Choose \(u\) so that
\(u \exp(-g^2u/2)\le1/2\). With probability at least \(1/2\) all
selected vertices thus have indegree at least \(\alpha u\).
The sample then contains a directed triangle. The expected number
of sampled triangles is their original number times
\(\binom u3/\binom v3\), giving the stated bound.

Now take \(\lambda>1/(1-\alpha)=2000/1307\) and suppose that such
avoiding colorings exist for arbitrarily large \(r\). Let \(D\) be
the oriented graph of arcs of \(M\) whose reverses are absent. Let \(b(x)\) be the number of
uncovered pairs at \(x\). The other vertices either receive an arc
from \(x\) in \(M\), send an unopposed arc to \(x\), or form an
uncovered pair with \(x\). By \eqref{pc:eq:outgoing}, the indegree
in \(D\) satisfies \(d_D^-(x)=N-1-d_M^+(x)-b(x)\ge N-1-h\binom{r}{t}/\binom{s}{t}-b(x)\).

Put \(g=(1-1/\lambda-\alpha)/4>0\). If at most \(\xi N^2\) pairs
are uncovered, deleting all vertices with \(b(x)>\sqrt\xi N\)
deletes at most \(2\sqrt\xi N\) vertices. The remaining graph has
minimum indegree at least \(N-1-h\binom{r}{t}/\binom{s}{t}-3\sqrt\xi N\). Choose \(\xi>0\)
small enough that its order is at least \(N/2\) and this bound is
at least \((\alpha+g)\) times its order for large \(N\). It then has
at least \(c_g N^3/8\) directed triangles.

The three directions on a directed triangle of \(D\) are the only
directions of \(M\) on those pairs. Such a triple contains neither
a bidirectional pair nor an in-star. By Claim~\ref{cl:agreement},
choose \(\rho\) sufficiently small and then \(r\) sufficiently large
so that there are at most \(\xi N^2\) uncovered pairs and fewer than
\(c_gN^3/8\) triples of this kind. This is a contradiction. Letting
\(\lambda\downarrow2000/1307\) proves the set-coloring upper bound.

Finally let \(k\to\infty\) and take the least \(r\) with \(\binom r2\ge k\).
Then \(\binom r2=k+O(\sqrt k)\). Encode the \(k\) ordinary colors by
distinct two-element labels. Two colors common to every label of a
copy force all those labels to be equal, so
\(R_3(P_n^{\mathrm{lin}};k,1,1)\le R_3(P_n^{\mathrm{lin}};r,2,2)\).
Applying the bound with \(s=t=2\) gives the ordinary upper estimate.
\end{theoremproof}

\section*{Funding}
This research did not receive any specific grant from funding agencies
in the public, commercial, or not-for-profit sectors.

\section*{Declaration of competing interest}
The authors declare that they have no competing interests.

\section*{Declaration of generative AI and AI-assisted technologies in the manuscript preparation process}
The authors used AI to assist with the organization and
wording of the manuscript and with checking and expanding steps in the
proofs. The authors are responsible for the mathematical arguments and
the final text.

\end{document}